\documentclass[12pt]{article}
\usepackage{graphicx}
\usepackage{graphics}
\usepackage{amsthm}
\usepackage{amssymb, amsmath}
\title{Mellin transforms with only critical zeros:  generalized Hermite functions}
\author{Mark W. Coffey\\
Department of Physics\\
Colorado School of Mines\\
Golden, CO  80401\\
USA\\
mcoffey@mines.edu}
\date{August 23, 2013}
\begin{document}
\maketitle
\baselineskip=25 pt
\begin{abstract}

We consider the Mellin transforms of certain generalized Hermite functions based upon certain generalized Hermite polynomials, characterized by a parameter $\mu>-1/2$.  We show that the transforms have polynomial factors whose zeros lie all on the critical line.  The polynomials with zeros only on the critical line are identified in terms of certain $_2F_1(2)$ 
hypergeometric functions, being certain scaled and shifted Meixner-Pollaczek polynomials. 
Other results of special function theory are presented. 

\end{abstract}
 
\vspace{.25cm}
\baselineskip=15pt
\centerline{\bf Key words and phrases}
\medskip 
Mellin transformation, generalized Hermite polynomials, hypergeometric function, critical line, zeros, functional equation, difference equation, reciprocity

\bigskip
\noindent
{\bf 2010 MSC numbers}
\newline{33C20, 33C45, 42C05, 44A20, 30D05}  

\baselineskip=25pt

\pagebreak
\centerline{\bf Introduction and statement of results}

\bigskip

In a series of papers, we are considering certain Mellin transforms comprised of classical orthogonal polynomials that yield polynomial factors with zeros only on the critical line Re $s=1/2$ or else only on the real axis.  Such polynomials have many important applications to
analytic number theory, in a sense extending the Riemann hypothesis.  For example,
using the Mellin transforms of Hermite functions, Hermite polynomials multiplied by a
Gaussian factor, Bump and Ng \cite{bumpng} were able to generalize
Riemann's second proof of the functional equation of the zeta function $\zeta(s)$, and to
obtain a new representation for it.  The polynomial factors turn out to be certain 
$_2F_1(2)$ Gauss hypergeometric functions \cite{coffeymellin}.  

In a different setting, the polynomials $p_n(x)= ~_2F_1(-n,-x;1;2)=(-1)^n ~_2F_1(-n,x+1;1;2)$ and $q_n(x)=i^n n! p_n(-1/2-ix/2)$ were studied \cite{kirsch}, and they directly correspond to the Bump and Ng polynomials with $s=-x$.  Kirschenhofer, Peth\"{o}, and Tichy considered combinatorial properties of $p_n$, and developed diophantine properties of them.  Their analytic results for $p_n$ include univariate and bivariate generating functions, and that its
zeros are simple, lie on the line $x=-1/2+it$, $t \in \mathbb{R}$, and that its zeros
interlace with those of $p_{n+1}$ on this line.  We may observe that these polynomials
may as well be written as $p_n(x)={{n+x} \choose n} ~_2F_1(-n,-x;-n-x;-1)$, or
$$p_n(x)={{(-1)^n 2^n \Gamma(n-x)} \over {n!\Gamma(-x)}} ~_2F_1\left(-n,-n;x+1-n;{1 \over 2}
\right),$$
where $\Gamma$ is the Gamma function.  In fact, combinatorial, geometrical, and coding aspects of $p_n(x)$ at integer argument had been noted in \cite{golomb} and \cite{stanton},
and Lemmas 2.2 and 2.3 of \cite{kirsch} correspond very closely to Lemmas 2 and 3,
respectively, of \cite{stanton}.

The Hermite polynomials being certain cases proportional to Laguerre polynomials $x^\delta L_n^{\pm 1/2}(x^2)$, $\delta=0$ or $1$, the generalization to Mellin transforms of Laguerre functions has been made \cite{bumpchoi,coffeymellin} and now the polynomial factors are a family of other $_2F_1(2)$ functions.  The Laguerre functions are ${\cal L}_n^\alpha(x) =x^{\alpha/2}e^{-x/2}L_n^\alpha(x)$, for $\alpha>-1$, and their Mellin transform is of the
form $M_n^\alpha(s)=2^{s+\alpha/2}\Gamma(s+\alpha/2)P_n^\alpha(s)$.  Mixed recursion
relations are known for the polynomials $P_n^\alpha$, as well as a generating function, and
they satisfy the functional equation $P_n^\alpha(s)=(-1)^n P_n^\alpha(1-s)$.

The generalized Mellin transform of Legendre, associated Legendre, Chebyshev, and
Gegenbauer functions has been investigated very recently elsewhere \cite{coffeylettunpub}.
These integral transforms are on the interval $[0,1]$ and lead to families of polynomials,
being certain terminating $_3F_2(1)$ functions, with zeros only on the critical line.

In this article, we study the Mellin transforms of certain generalized Hermite functions, and
are able to identify the resulting polynomial factors in terms of certain hypergeometric functions $_2F_1(2)$.  The key result is that these transforms possess zeros only on the critical line. 

We use standard notation.
Let $_pF_q$ be the generalized hypergeometric function, $(a)_n=\Gamma(a+n)/\Gamma(a)=(-1)^n {{\Gamma(1-a)} \over {\Gamma(1-a-n)}}$ the Pochhammer symbol, 
and $B(x,y)=\Gamma(x)\Gamma(y)/\Gamma(x+y)$ the Beta function.  



Szeg\"{o} in a problem introduced an orthogonal polynomial sequence $\{H_n^\mu(x)\}$
\cite{szego} (p. 377).  These polynomials were also studied early on by Chihara \cite{chihara}.
They may be written as
$$H_{2n}^\mu(x)=(-1)^n 2^{2n} n! L_n^{\mu-1/2}(x^2), \eqno(1.1a)$$
and
$$H_{2n+1}^\mu(x)=(-1)^n 2^{2n+1} n! xL_n^{\mu+1/2}(x^2), \eqno(1.1b)$$
where $\mu>-1/2$ and $L_n^\alpha$ are Laguerre polynomials.  
In turn, the Laguerre polynomials are expressible in terms of the confluent hypergeometric
function $_1F_1$ as 
$$L_n^\alpha(x) ={{n+\alpha} \choose n} ~_1F_1(-n;\alpha+1;x). \eqno(1.2)$$
The polynomials $H_n^\mu(x)$ have been written by other authors with a different normalization. 
In particular, Rosenblum \cite{rosenblum} does not include the $2^p$ factors.  
The polynomials are othogonal with respect to the weight function $|x|^{2\mu}e^{-x^2}$, for
$x \in (-\infty,\infty)$, such that
$$\int_{-\infty}^\infty H_m^\mu(x)H_n^\mu(x)|x|^{2\mu}e^{-x^2}dx=2^{2n}\left[{n\over 2}\right]
\Gamma\left(\left[{{n+1}\over 2}\right]+\mu+{1 \over 2}\right)\delta_{mn},$$
with $[x]$ the greatest integer function.

Just as the Hermite polynomials are connected with the excited state wavefunctions and
algebraic properties of the quantum mechanical simple harmonic oscillator, the generalized
Hermite polynomials may be used to develop the calculus of the Bose-like oscillator
\cite{rosenblum}.  Generalized Hermite functions may also be shown to be the 
eigenfunctions of a generalized Fourier transform, with eigenvalues $\pm 1$ and $\pm i$.

Finite sums for Laguerre polynomials and Kummer-type relations for the hypergeometric function
$_2F_2$ are established in \cite{coffeysums}.  Convolution sums and other results for Laguerre
polynomials are given in \cite{coffeyjpa}.  As we briefly show in the third section of the
paper, these relations may be usefully transcribed to the generalized Hermite polynomials.

The polynomials $H_n^\mu(x)$ satisfy the recurrence relation
$$H_{n+1}^\mu(x)=2xH_n^\mu(x)-2(n+\theta_n)H_{n-1}^\mu(x), ~~~~~~n \geq 0, \eqno(1.3)$$
as well as the differential equation
$$xy''+2(\mu-x^2)y'+(2nx-\theta_n x^{-1})y=0, \eqno(1.4)$$
where $\theta_{2m}=0$ and $\theta_{2m+1}=2\mu$.  

Using a generating function of Laguerre polynomials, one for the generalized Hermite 
polynomials may be developed \cite{chihara},
$$(1+2xw+4w^2)(1+4w^2)^{-\mu+3/2}\exp[4x^2w^2(1+4w^2)^{-1}]
= \sum_{n=0}^\infty H_n^\mu(x) {w^n \over {[n/2]!}}.  \eqno(1.5)$$
The polynomials $H_n^\mu(x)$ have application in Gauss-generalized-Hermite quadrature \cite{shao}.


Here we consider Mellin transformations 
$$({\cal M}f)(s)=\int_0^\infty f(x)x^s {{dx} \over x}.$$
For properties of the Mellin transform, we mention \cite{butzer}.

We put, for Re $s>0$,
$$M_n^\mu(s) \equiv \int_0^\infty x^{s-1} H_n^\mu(x) x^\mu e^{-x^2/2}dx.  \eqno(1.6)$$
We let $p_n^\mu(s)$ denote the polynomial factor of $M_n^\mu(s)$.
Clearly (1.6) is the $\mu \neq 0$ extension of the Mellin transform of Hermite functions.

{\bf Proposition 1}.  The Mellin transforms (1.3) are given by
$$M_{2n}^\mu(s)=(-1)^n2^{2n}2^{(\mu+s)/2-1}\left(\mu+{1 \over 2}\right)_n \Gamma\left({{\mu+s} \over 2}\right)~_2F_1\left(-n,{{\mu+s} \over 2};\mu+{1 \over 2};2\right), \eqno(1.7a)$$
and
$$M_{2n+1}^\mu(s)=(-1)^n2^{2n}2^{(\mu+s+1)/2}\left(\mu+{3 \over 2}\right)_n \Gamma\left({{\mu+s+1} \over 2}\right)~_2F_1\left(-n,{{\mu+s+1} \over 2};\mu+{3 \over 2};2\right), \eqno(1.7b)$$

{\bf Proposition 2}.  The Mellin transforms (1.3) satisfy the recursions
$$M_{2m+1}^\mu(s)=2M_{2m}^\mu(s+1)-4m M_{2m-1}^\mu(s), \eqno(1.8a)$$
and
$$M_{2m+2}^\mu(s)=2M_{2m+1}^\mu(s+1)-2(2m+2\mu+1)m M_{2m}^\mu(s). \eqno(1.8b)$$

{\bf Proposition 3}. A generating function for the transforms (1.6) is
$$\sum_{n=0}^\infty M_n^\mu(s) {t^n \over {[n/2]!}}=2^{(\mu+s-3)/2} {{(1+4t^2)^{s/2-1}} \over
{(1-4t^2)^{(\mu+s+1)/2}}} \left[\sqrt{2(1-16t^4)} \Gamma\left({{\mu+s} \over 2}\right)
+4t \Gamma\left({{\mu+s+1} \over 2}\right)\right].$$

{\bf Proposition 4}.  (Reciprocity) For positive integers $m$ and $n$,
$$\left(\mu+{1 \over 2}\right)_{m}p_{2n}^\mu(-2m-\mu)=\left(\mu+{1 \over 2}\right)_{n}p_{2m}^\mu(-2n-\mu),$$
and 
$$\left(\mu+{3 \over 2}\right)_{m}p_{2n+1}^\mu(-2m-1-\mu)=\left(\mu+{3 \over 2}\right)_{n}p_{2m+1}^\mu(-2n-1-\mu).$$

{\bf Theorem 1}.
The polynomials $p_n^\mu(s)$, of degree $n$, satisfy the functional equation
$p_n^\mu(s)=(-1)^n p_n^\mu(1-s)$.  These polynomials have zeros only on the
critical line.  Further, all zeros $\neq 1/2$ occur in complex conjugate pairs.


The following section of the paper contains the proof of these Propositions. 
After that are presented various special function theory results and some discussion.

\medskip
\centerline{\bf Proof of Propositions and of Theorem 1}
\medskip

{\it Proposition 1}.  Using the definition (1.1) for $H_{2n}^\mu$ and the series form of (1.2)  
for the Laguerre polynomials in (1.6) gives
$$M_{2n}^\mu(s)=(-1)^n 2^{2n}n! \int_0^\infty x^{s+\mu-1} e^{-x^2/2}L_n^{\mu-1/2}(x^2)dx$$
$$=(-1)^n2^{2n}\left(\mu+{1 \over 2}\right)_n\sum_{j=0}^n {{(-n)_j} \over {(\mu+1/2)_j}}{1 \over
{j!}} \int_0^\infty x^{s+\mu+2j-1} e^{-x^2/2}dx$$
$$=(-1)^n2^{2n}\left(\mu+{1 \over 2}\right)_n\sum_{j=0}^n {{(-n)_j} \over {(\mu+1/2)_j}}{1 \over
{j!}} \Gamma\left(j+{{\mu+s} \over 2} \right)2^{(\mu+s)/2+j-1}$$
$$=(-1)^n2^{2n}\left(\mu+{1 \over 2}\right)_n2^{(\mu+s)/2-1}\Gamma\left({{\mu+s} \over 2} \right)
\sum_{j=0}^n {{(-n)_j} \over {(\mu+1/2)_j}}\left({{\mu+s} \over 2}\right)_j{2^j \over {j!}}$$
$$=(-1)^n2^{2n}2^{(\mu+s)/2-1}\left(\mu+{1 \over 2}\right)_n \Gamma\left({{\mu+s} \over 2}\right)~_2F_1\left(-n,{{\mu+s} \over 2};\mu+{1 \over 2};2\right).$$
Part (b) proceeds very similarly.  \qed

{\it Proposition 2}.  The recursions follow directly from (1.1) and (1.3).  \qed

{\it Proposition 3}.  The generating function follows by using (1.5) in the definition (1.6).  \qed

{\it Remark}.  Expanding the generating function so obtained in powers of $t$ is another way in
which to find the explicit transform expressions of Proposition 1.

{\it Proposition 4}.  This follows from the symmetry of the $_2F_1$ function in its two
numerator parameters.  \qed

{\it Remark}.  We highly expect that the reciprocity relation and functional equation for
$p_n^\mu(s)$ are connected with the properties of an Ehrhart polynomial, possibly shifted and/or
scaled, corresponding to some polytope.

{\it Theorem 1}.  By Proposition 1, the polynomials $p_n^\mu(s)$ are of degree $n$ and have
real coefficients for real values of $\mu$.  The functional equation follows immediately by
transforming the $_2F_1(2)$ function.  In particular we apply (e.g., \cite{grad} p. 1043)
$$_2F_1(\alpha,\beta;\gamma;z)=(1-z)^{-\alpha} ~_2F_1\left(\alpha,\gamma-\beta;\gamma;{z \over {z-1}}\right).$$

That the polynomial zeros occur only on Re $s=1/2$ will follow from the following 
difference equations, which we next derive.  For $n$ even,
$$[2(n+\mu)+1](s+\mu-2)p_n^\mu(s)-(s+\mu)(s+\mu-2)p_n^\mu(s+2)+[(s-2)(s-1)+(1-\mu)\mu]
p_n^\mu(s-2)=0,$$
and for $n$ odd,
$$[2(n+\mu)-1](s+\mu-1)p_n^\mu(s)-(s+\mu+1)(s+\mu-1)p_n^\mu(s+2)+[(s-2)(s-1)-(1+\mu)\mu]
p_n^\mu(s-2)=0.$$
We put $H_n^\mu(x)=x^{-\mu}e^{x^2/2}f(x)$ into the differential equation (1.4).
Then for $n$ even there results
$$e^{x^2/2}x^{-\mu-1}\left\{[(1-\mu)\mu+x^2+2(n+\mu)x^2-x^4]f(x)+x^2f''(x)\right\}=0,$$
and
$$e^{x^2/2}x^{-\mu-1}\left\{-[(1+\mu)\mu-(1+2(n+\mu))x^2+x^4]f(x)+x^2f''(x)\right\}=0$$
for $n$ odd.
The quantity in curly brackets is zero, and multiplying it by $x^{s-1}$ and
integrating by parts twice and shifting $s\to s-2$, for $n$ even we find the difference
equation of the Mellin transforms,
$$[2(n+\mu)+1]M_n^\mu(s)-M_n^\mu(s+2)+[(1-\mu)\mu+(s-2)(s-1)]M_n^\mu(s-2)=0.$$
By using the $s$-dependent factors of the transforms (1.7a),
$M_n^\mu(s) \propto 2^{s/2}\Gamma\left({{s+\mu} \over 2}\right)p_n^\mu(s)$,
and the functional equation $\Gamma(z+1)=z\Gamma(z)$, we find the stated difference
equation for the polynomial factors when $n$ is even.  For $n$ odd, the steps are very
similar.

We then use shifted polynomials $q(s)=p_n^\mu(s+1/2)$, so that $p_n^\mu(s)=q(s-1/2)$.
Then, with a translation $s \to s+1/2$, we find
$$[2(n+\mu)+1]\left(s+\mu-{3 \over 2}\right)q(s)-\left(s+\mu+{1 \over 2}\right)\left(s+\mu-{3 \over 2}\right)q(s+2)$$
$$+\left[\left(s-{3 \over 2}\right)\left(s-{1 \over 2}\right)+(1-\mu)\mu \right]q(s-2)=0,$$
for $n$ even, and
$$[2(n+\mu)-1]\left(s+\mu-{1 \over 2}\right)q(s)-\left(s+\mu+{3 \over 2}\right)\left(s+\mu-{1 \over 2}\right)q(s+2)$$
$$+\left[\left(s-{3 \over 2}\right)\left(s-{1 \over 2}\right)-(1+\mu)\mu \right]q(s-2)=0,$$
for $n$ odd.  It follows that if $r_k$ is a root of $q$, $q(r_k)=0$, that
$$\left(r_k+\mu+{1 \over 2}\right)q(r_k+2)=\left(r_k-\mu-{1 \over 2}\right)q(r_k-2),$$
when $n$ is even, and
$$\left(r_k+\mu+{3 \over 2}\right)q(r_k+2)=\left(r_k-\mu-{3 \over 2}\right)q(r_k-2),$$
when $n$ is odd.  The equality of the absolute value of both sides of these equations
provides a necessary condition that Re $r_k=0$.  I.e., the zeros of $q(s)$ are pure
imaginary and thus the zeros of $p_n^\mu(s)$ are on the critical line.  \qed 

\pagebreak
\centerline{\bf Other results}
\medskip

Let $\varepsilon(m)=0$  when $m$ is even and $=1$ when $m$ is odd.  Then we may compactly
write from (1.1)
$$H_m^\mu(x)=(-1)^{{(m-\varepsilon)}/2}2^m\left({{m-\varepsilon} \over 2}\right)! x^\varepsilon
L_{{m-\varepsilon} \over 2}^{\mu-1/2+\varepsilon}(x^2).  \eqno(2.1)$$
We let $J_n$ be the Bessel function of the first kind of order $n$.
We then have the following integral representation.
{\newline \bf Proposition 5}.  For $\mu>-1/2$,
$$H_m^\mu(x)=(-1)^{{(m-\varepsilon)}/2}2^m e^{x^2} x^{1/2-\mu} \int_0^\infty t^{m/2+\mu/2-1/4}
J_{\mu-1/2+\varepsilon}(2x\sqrt{t})e^{-t}dt. \eqno(2.2)$$

{\it Proof}.  This follows from (2.1) and the representation (\cite{andrews}, p. 286) for $\alpha>-1$
$$L_n^\alpha(x)={1 \over {n!}}e^x x^{-\alpha/2} \int_0^\infty t^{n+\alpha/2}J_\alpha(2\sqrt{xt})
e^{-t}dt.$$
\qed

{\it Remark}.  As the representation (2.2) is conditionally convergent, it does not appear
to be directly useful for application to the transforms (1.6).

{\bf Proposition 6}.   For $\beta>0$ and $\mu>-1/2$,
$$H_m^{\mu+\beta}(x)={{\Gamma\left({{m+\varepsilon+1} \over 2}+\mu+\beta\right)} \over {\Gamma(\beta)\Gamma\left(
{{m+\varepsilon+1} \over 2}+\mu\right)}} \int_0^1 t^{\mu-1/2+\varepsilon/2}(1-t)^{\beta-1} 
H_m^\mu(x \sqrt{t})dt$$
$$=2{{\Gamma\left({{m+\varepsilon+1} \over 2}+\mu+\beta\right)} \over {\Gamma(\beta)\Gamma\left(
{{m+\varepsilon+1} \over 2}+\mu\right)}}\int_0^1 u^{2\mu+\varepsilon}(1-u^2)^{\beta-1} H_m^\mu
(xu)du.$$
Here $\varepsilon=0$ when $m$ is even and $=1$ when $m$ is odd.

{\it Proof}.  This follows from (1.1) and Koshlyakov's formula (\cite{kosh}, \cite{lebedev}, 
p. 94)
$$L_n^{\alpha+\beta}(x)={{\Gamma(n+\alpha+\beta+1)} \over {\Gamma(\beta)\Gamma(n+\alpha+1)}}
\int_0^1 t^\alpha (1-t)^{\beta-1} L_n^\alpha(xt)dt, ~~~~\alpha >-1, ~~~~\beta>0.$$
\qed

{\it Remark}.  This result contains the very special case of $\mu=2n+\epsilon=0$ and $\beta \to \alpha+1/2$ in 7.372 in \cite{grad}, p. 836.  Indeed this corresponds to Uspensky's formula
(e.g., \cite{lebedev}, p. 94). 

{\bf Proposition 7}.  For $\mu>-1/2$,
$$\sum_{m=0}^n {{(-1)^m} \over {2^{2m}m!}}H_{2m}^{\mu+1/2}(x)= {{(-1)^n} \over {2^{2n}n!}}H_{2n}^{\mu+3/2}(x)$$
and
$$\sum_{m=0}^n {{(-1)^m} \over {2^{2m}m!}}H_{2m+1}^{\mu-1/2}(x)= {{(-1)^n} \over {2^{2n}n!}}H_{2n+1}^{\mu+1/2}(x).$$

{\it Proof}. This follows from (1.1) and the relation (e.g., \cite{grad}, p. 1038)
$\sum_{m=0}^n L_m^\alpha(x)=L_n^{\alpha+1}(x)$. \qed

We may write the generalized Hermite polynomials in terms of the ordinary Hermite polynomials,
and vice versa, as described next.
{\newline \bf Proposition 8}.  For $\mu>-1/2$,
(a) 
$$\sum_{j=0}^n {n \choose j}(-1)^j4^j (-\mu)_j H_{2(n-j)}^\mu(x)=H_{2n}(x),$$
(b)
$$\sum_{j=0}^n {n \choose j}(-1)^j4^j (\mu)_j H_{2(n-j)}(x)=H_{2n}^\mu(x),$$
(c)
$$\sum_{j=0}^n {n \choose j}(-1)^j4^j (-\mu)_j H_{2(n-j)+1}^\mu(x)=H_{2n+1}(x),$$
(d)
$$\sum_{j=0}^n {n \choose j}(-1)^j4^j (\mu)_j H_{2(n-j)+1}(x)=H_{2n+1}^\mu(x),$$
thus
(e)
$$H_n^\mu(x)=\sum_{j=0}^{[n/2]} {{[n/2]} \choose j}(-1)^j 4^j (\mu)_j H_{n-2j}(x)$$
and
$$H_n(x)=\sum_{j=0}^{[n/2]} {{[n/2]} \choose j}(-1)^j 4^j (-\mu)_j H_{n-2j}^\mu(x).$$

{\it Proof}.  These relations follows from (1.1) and the property
$$L_n^\alpha(x)=\sum_{j=0}^n {{(\alpha-\beta)_j} \over {j!}}L_{n-j}^\beta(x).$$
Thus parts (a) and (b) follow from
$$\sum_{j=0}^n (\alpha-\mu)_j {n \choose j} (-1)^j 4^j H_{2(n-j)}^\mu(x)=H_{2n}^\alpha(x),$$
and (c) and (d) from
$$\sum_{j=0}^n (\alpha-\mu)_j {n \choose j} (-1)^j 4^j H_{2(n-j)+1}^\mu(x)=H_{2n+1}^\alpha(x).$$
\qed

Proposition 8 implies a summation relation for certain $_2F_1(2)$ functions.  Stated otherwise, we have the following, with $M_n(s)=M_n^0(s)$.
{\newline \bf Corollary}.  For $\mu>-1/2$,
$$M_n^\mu(s)=\sum_{j=0}^{[n/2]} {{[n/2]} \choose j}(-1)^j 4^j (\mu)_j M_{n-2j}(s).$$

We have the following convolution sums for generalized Hermite polynomials.
{\bf Proposition 9}.  (a)
$$H_{2n}^{\alpha+\beta+1/2}(\sqrt{x^2+y^2})=\sum_{k=0}^n {n \choose k}H_{2(n-k)}^\alpha(x)
H_{2k}^\beta(y),$$
(b)
$${{xy} \over \sqrt{x^2+y^2}}H_{2n+1}^{\alpha+\beta+3/2}(\sqrt{x^2+y^2})={1 \over 2}\sum_{k=0}^n {n \choose k}H_{2(n-k)+1}^\alpha(x)H_{2k+1}^\beta(y),$$
(c) 
$$\sum_{j=1}^{n-1}{{n-1} \choose {j-1}}H_{2(j-1)}^\alpha(x)H_{2(n-j-1)}^{k+3/2}(y)$$
$$=-{1 \over 4}{{k!} \over y^{2(k+1)}}\left[\sum_{\ell_1=0}^k {y^{2\ell_1} \over {\ell_1!}}
H_{2(n-1)}^\alpha(x)-\sum_{\ell_2=0}^k {y^{2\ell_2} \over {\ell_2!}}
H_{2(n-1)}^{\alpha+\ell_2}(\sqrt{x^2+y^2})\right],$$
(d)
$${2 \over {xy}}\sum_{j=1}^{n-1}{{n-1} \choose {j-1}}H_{2j-1}^\alpha(x) H_{2(n-j)-1}^{k+1/2}(y)$$
$$=-{{k!} \over y^{2(k+1)}}\left[{1 \over x}\sum_{\ell_1=0}^k {y^{2\ell_1} \over {\ell_1!}}
H_{2n-1}^\alpha(x)-{1 \over \sqrt{x^2+y^2}}\sum_{\ell_2=0}^k {y^{2\ell_2} \over {\ell_2!}}
H_{2n-1}^{\alpha+\ell_2}(\sqrt{x^2+y^2})\right],$$
(e)
$$H_{2n}^{\alpha_1+\alpha_2+\ldots+\alpha_k-1/2}\left(\sqrt{x_1^2+x_2^2+\ldots+x_k^2}\right)$$
$$=n!\sum_{i_1+i_2+\ldots+i_k=n} {{H_{2i_1}^{\alpha_1+1/2}(x_1)H_{2i_2}^{\alpha_2+1/2}(x_2)\cdots
H_{2i_k}^{\alpha_k+1/2}(x_k)} \over {i_1!i_2!\cdots i_k!}},$$
and (f)
$${{x_1x_2\cdots x_k} \over \sqrt{x_1^2+x_2^2+\ldots+x_k^2}} H_{2n+1}^{\alpha_1+\alpha_2+\ldots+\alpha_k-3/2}\left(\sqrt{x_1^2+x_2^2+\ldots+x_k^2}\right)$$
$$={{n!} \over 2^{k-1}}\sum_{i_1+i_2+\ldots+i_k=n} {{H_{2i_1+1}^{\alpha_1-1/2}(x_1)H_{2i_2+1}^{\alpha_2-1/2}(x_2)\cdots
H_{2i_k+1}^{\alpha_k-1/2}(x_k)} \over {i_1!i_2!\cdots i_k!}}.$$
In (c) and (d), $k>0$ is an integer.

{\it Proof}.  (a) and (b) follow from (1.1) and the well known relation (e.g., \cite{grad}, p. 1038)
$$L_n^{\alpha+\beta+1}(x+y)=\sum_{k=0}^n L_{n-k}^\alpha(x)L_k^\beta(y).$$

(c) and (d) follow from (1.1) and Proposition 3 of \cite{coffeyjpa}:
For $\alpha >-1$ and $k$ a positive integer, 
$$\sum_{j=1}^{n-1} {{L_{j-1}^\alpha(x)L_{n-j-1}^{k+1}(y)} \over {n-j}}
={{k!} \over y^{k+1}}\left[\sum_{\ell_1=0}^k{y^{\ell_1} \over \ell_1!}
L_{n-1}^\alpha(x)-\sum_{\ell_2=0}^k {y^{\ell_2} \over \ell_2!}
L_{n-1}^{\alpha+\ell_2}(x+y)\right].$$
(e) and (f) follow from (1.1) and (e.g., \cite{grad}, p. 1039)
$$L_n^{\alpha_1+\alpha_2+\ldots+\alpha_k+k-1}(x_1+x_2+\ldots+x_k)
=\sum_{i_1+i_2+\ldots+i_k=n} L_{i_1}^{\alpha_1}(x_1)L_{i_2}^{\alpha_2}(x_2)\ldots L_{i_k}^{\alpha_k}(x_k).$$
\qed

{\bf Proposition 10}.   
(a)  For $|t|<1$,
$$\sum_{m=0}^\infty {{(-1)^m} \over {4^m m!}}H_{2m}^{\alpha-m+1/2}(x)t^m=(1+t)^\alpha e^{-x^2t},$$
$$\sum_{m=0}^\infty {{(-1)^m} \over {4^m m!}}H_{2m+1}^{\alpha-m-1/2}(x)t^m=2x(1+t)^\alpha e^{-x^2t},$$
(b)
$$H_{2n}^\alpha(\sqrt{x^2+y^2})=e^{y^2}\sum_{k=0}^\infty {{(-1)^k} \over {k!}}y^{2k} H_{2n}^{\alpha+k}(x),$$
$$H_{2n+1}^\alpha(\sqrt{x^2+y^2})=e^{y^2}\sum_{k=0}^\infty {{(-1)^k} \over {k!}}y^{2k} H_{2n+1}^{\alpha+k}(x),$$
(c) 
$$\sum_{m=0}^\infty {t^m \over {(\alpha+1)_m}} {{(-1)^m} \over {4^m m!}}H_{2m}^{\alpha+1/2}(x)
=\Gamma(\alpha+1){e^t \over {t^{\alpha/2} x^\alpha}}J_\alpha(2x\sqrt{t}),$$
$$\sum_{m=0}^\infty {t^m \over {(\alpha+1)_m}} {{(-1)^m} \over {2^{2m+1} m!}}H_{2m+1}^{\alpha+1/2}(x)
=\Gamma(\alpha+1){e^t \over {t^{\alpha/2} x^{\alpha-1}}}J_\alpha(2x\sqrt{t}),$$
(d) for $x>0$ and $\gamma$ the Euler constant,
$$\sum_{n=1}^\infty {{(-1)^n} \over {4^n n n!}}H_{2n}^{1/2}(x)=-2\ln x -\gamma,$$
$$\sum_{n=1}^\infty {{(-1)^n} \over {2^{2n+1} n n!}}H_{2n+1}^{-1/2}(x)=-x(2\ln x +\gamma),$$
(e)
$${{(-1)^m } \over {4^m m!}}H_{2m}^{\beta+1/2}(\tau x)=\sum_{n=0}^{m} {{(-1)^n } \over {4^n n!}}
{{\beta+m} \choose {\beta-n}}\tau^{2n}(1-\tau^2)^{m-n} H_{2n}^{\beta+1/2}(x),$$
and (f)
$${{(-1)^m } \over {2^{2m+1} m!}}{1 \over \tau}H_{2m+1}^{\beta-1/2}(\tau x)=\sum_{n=0}^m
{{(-1)^n } \over {2^{2n+1} n!}}{{\beta+m} \choose {\beta-n}}\tau^{2n}(1-\tau^2)^{m-n} H_{2n+1}^{\beta-1/2}(x).$$
(e) and (f) continue to hold for $\tau \to 0$, when only the $n=0$ term remains on the right side.  

{\it Proof}.  (a) follows from (1.1) and the generating function
$$\sum_{m=0}^\infty L_m^{\alpha-m}(x)t^m=(1+t)^\alpha e^{-xt}.$$
For (b),
$$L_n^\alpha(x+y)=e^y \sum_{k=0}^\infty {{(-1)^k} \over {k!}}y^k L_n^{\alpha+k}(x)$$
is used.  (c) follows from
$$\sum_{m=0}^\infty {t^m \over {(\alpha+1)_m}} L_m^\alpha(x)=e^t ~_0F_1(-;\alpha+1;-xt)
=\Gamma(\alpha+1){e^t \over {(xt)^{\alpha/2}}}J_\alpha(2\sqrt{xt}),$$
where $_0F_1$ is a confluent hypergeometric function with a single denominator parameter.
(d) is based upon
$$\sum_{n=1}^\infty {{L_n(x)} \over n}=-\ln x-\gamma, ~~~~~~x>0.$$ 
(e) and (f) use (1.1) and
$$L_m^\beta(\tau x)=\sum_{n=0}^m {{\beta+m} \choose {\beta-n}}\tau^n (1-\tau)^{m-n} L_n^\beta(x).$$
\qed

The ordinary Hermite polynomials have the alternative representation 
$$H_n(x)=(2x)^n ~_2F_0\left(-[n/2],-[n/2]+(-1)^n/2;-;-{1 \over x^2}\right),$$
so that another way to write the generalized Hermite polynomials is
$$H_n^\mu(x)=(2x)^n ~_2F_0\left(-[n/2],-[n/2]-\mu+(-1)^n/2;-;-{1 \over x^2}\right).$$
In turn, the Mellin transforms (1.6) evaluate in terms of $_2F_1(1/2)$ functions.
These functions are to be expected, as transformations of the hypergeometric factors in (1.7),
$$_2F_1\left(-n,{{\mu+s} \over 2};\mu+{1 \over 2};2\right)=(-1)^n2^n 
{{\Gamma\left(\mu+{1 \over 2}\right)\Gamma\left({{\mu+s} \over 2}+n\right)} \over
{\Gamma\left({{\mu+s} \over 2}\right)\Gamma\left(\mu+{1 \over 2}+n\right)}}.$$
$$\times ~_2F_1\left(-n,{1 \over 2}-n-\mu;1-n-{{(\mu+s)} \over 2};{1 \over 2}\right).$$ 
$$_2F_1\left(-n,{{\mu+s+1} \over 2};\mu+{3 \over 2};2\right)=(-1)^n2^n 
{{\Gamma\left(\mu+{3 \over 2}\right)\Gamma\left({{\mu+s+1} \over 2}+n\right)} \over
{\Gamma\left({{\mu+s+1} \over 2}\right)\Gamma\left(\mu+{3 \over 2}+n\right)}}.$$
$$\times ~_2F_1\left(-n,-{1 \over 2}-n-\mu;1-n-{{(\mu+s+1)} \over 2};{1 \over 2}\right).$$

The Meixner-Pollaczek polynomials $P_n^\lambda$ have the hypergeometric form (e.g., \cite{andrews}, p. 348)
$$P_n^\lambda(x;\phi)={{(2\lambda)_n} \over {n!}} e^{in\phi} ~_2F_1\left(-n,\lambda+ix;2\lambda;
1-e^{-2i\phi}\right), ~~~~\lambda>0, ~~~~0<\phi<\pi.$$
Therefore, comparing with (1.7), the polynomial factors there are proportional to
$$P_n^{(\mu+\varepsilon)/2+1/4}\left[i\left({1 \over 4}-{s \over 2}\right);{\pi \over 2}\right]
={{(\mu+1/2+\varepsilon)_n} \over {n!}} i^n ~_2F_1\left(-n,{{\mu+s+\varepsilon} \over 2};
\mu+{1 \over 2}+\varepsilon;2\right).$$
This identification provides an alternative method to prove Theorem 1, and, in fact, to
show that the polynomial zeros on the critical line are also simple.
First, the Meixner-Pollaczek polynomials are orthogonal with respect to the weight 
$$\left|\Gamma\left(\lambda+ix\right)\right|^2e^{(2\phi-\pi)x}.$$
The Plancherel relation for Mellin transforms provides none other than this weight,
$$\left|\Gamma\left({{\mu+s+\varepsilon} \over 2}\right)\right|^2.$$
Then the polynomial factors $p_n^\mu(1/2+it)$ form an orthogonal family with respect to the
corresponding measure.  By standard results on orthogonal polynomials \cite{szego}, the zeros 
are simple, and the zeros of $p_n^\mu(1/2+it)$ and $p_{n+1}^\mu(1/2+it)$ separate each other.  



\pagebreak


\begin{thebibliography}{99}
\bibitem{nbs}M. Abramowitz and I. A. Stegun,
{Handbook of Mathematical Functions, Washington, National Bureau of Standards (1964).}
\bibitem{andrews}G. E. Andrews, R. Askey, and R. Roy, 
{Special Functions, Cambridge University Press (1999).}
\bibitem{bumpchoi}D. Bump, K.-K. Choi, P. Kurlberg, and J. Vaaler,
{A local Riemann hypothesis, I, Math. Z. {\bf 233}, 1-19 (2000).}
\bibitem{bumpng}D. Bump and E. K.-S. Ng,
{On Riemann's zeta function, Math. Z. {\bf 192}, 195-204 (1986).}
\bibitem{butzer}P. Butzer and S. Jansche,
{A direct approach to the Mellin transform, J. Fourier Analysis Appls. {\bf 3}, 325-376 
(1997).}
\bibitem{chihara}T. Chihara,
{An introduction to orthogonal polynomials, Gordon and Breach (1978).}
\bibitem{coffeymellin}M. W. Coffey,
{Special functions and the Mellin transforms of Laguerre and Hermite functions, 
Analysis {\bf 27}, 95-108 (2007).}  
\bibitem{coffeyxi}M. W. Coffey,
{Theta and Riemann xi function representations from harmonic oscillator eigenfunctions,
Phys. Lett. A {\bf 362}, 352-356 (2007).}
\bibitem{coffeylettunpub}M. W. Coffey and M. C. Lettington,
{unpublished (2013).}
\bibitem{coffeysums}M. W. Coffey
{On finite sums of Laguerre polynomials, Rocky Mtn. J. Math. {\bf 41}, 79-93 (2011).}
\bibitem{coffeyjpa}M. W. Coffey,
{Properties and possibilities of quantum shapelets, J. Phys. A {\bf 39}, 877-887 (2006).}
\bibitem{golomb}S. W. Golomb and L. R. Welch,
{Perfect codes in the Lee metric and the packing of polyominoes, SIAM J. Appl. Math.
{\bf 18}, 302-317 (1970).}
\bibitem{grad}I. S. Gradshteyn and I. M. Ryzhik,
{Table of Integrals, Series, and Products, Academic Press, New York (1980).}
\bibitem{kirsch}P. Kirschenhofer, A. Peth\"{o}, and R. F. Tichy,
{On analytical and diophantine properties of a family of counting polynomials, 
Acta Sci. Math. (Szeged) {\bf 65}, 47-59 (1999).}
\bibitem{kosh}N. S. Koshlyakov,
{On Sonine's polynomials, Mess. Math. {\bf 55}, 152 (1926).}  
\bibitem{lebedev}N. N. Lebedev,
{Special functions and their applications, Dover (1972).}
\bibitem{rosenblum}M. Rosenblum,
{Generalized Hermite polynomials and the Bose-like oscillator calculus, Oper.
Theory Adv. Appl. {\bf 73}, 369-396 (1994).}
\bibitem{shao}T. S. Shao, T. C. Chen, and R. M. Frank,
{Tables of zeros and gaussian weights of certain associated Laguerre polynomials and the related
generalized Hermite polynomials, {\bf 18}, 598-616 (1964).}
\bibitem{stanton}R. G. Stanton and D. D. Cowan,
{Note on a ``square" functional equation, SIAM Rev. {\bf 12}, 277-279 (1970).}
\bibitem{szego}G. Szeg\"{o},
{Orthogonal Polynomials, Vol. 23 of AMS colloquium Publications, American Mathematical
Society, Providence, RI (1975).}
\end{thebibliography}
\end{document}